\documentclass[
a4paper,
11pt,
twoside,
]{article}
\usepackage[utf8]{inputenc}
\usepackage[T1]{fontenc}
\usepackage[osf,sc,noBBpl]{mathpazo} 
\usepackage{microtype}
\usepackage[english]{babel}
\usepackage{geometry}
\usepackage{amsmath,amsfonts,amssymb,amsthm}

\usepackage{etoolbox}

\usepackage{graphicx}
\usepackage[dvipsnames]{xcolor}

\theoremstyle{plain}
\newtheorem{theorem}{Theorem}[section]

\newtheorem{remark}[theorem]{Remark}

\usepackage{mathtools}
\usepackage{amssymb}
\usepackage{siunitx}
\usepackage{booktabs}

\newcommand{\dd}{\mathop{}\!\mathrm{d}}

\renewcommand{\restriction}{\raise-.5ex\hbox{\ensuremath{\upharpoonright}}}

\definecolor{webgreen}{rgb}{0,.5,0}
\definecolor{webbrown}{rgb}{.6,0,0}
\definecolor{myblue}{rgb}{0,0.25,0.5}

\usepackage[
  colorlinks=true,
  breaklinks=true,
  bookmarksnumbered,
  pdfhighlight=/O,
  urlcolor=webbrown,
  linkcolor=myblue,
  citecolor=webgreen,
  hyperfootnotes=true,
]{hyperref}
\newcommand{\email}[1]{\href{mailto:#1}{\texttt{#1}}}
\makeatletter
\newcommand{\pagerefstar}{\@pagerefstar}
\makeatother

\usepackage[capitalize,nameinlink,sort]{cleveref}[0.19]
\crefname{section}{section}{sections}
\crefname{subsection}{subsection}{subsections}
\Crefname{figure}{Figure}{Figures}
\crefformat{equation}{\textup{#2(#1)#3}}
\crefrangeformat{equation}{\textup{#3(#1)#4--#5(#2)#6}}
\crefmultiformat{equation}{\textup{#2(#1)#3}}{ and \textup{#2(#1)#3}}
{, \textup{#2(#1)#3}}{, and \textup{#2(#1)#3}}
\crefrangemultiformat{equation}{\textup{#3(#1)#4--#5(#2)#6}}%
{ and \textup{#3(#1)#4--#5(#2)#6}}{, \textup{#3(#1)#4--#5(#2)#6}}{, and \textup{#3(#1)#4--#5(#2)#6}}
\Crefformat{equation}{#2Equation~\textup{(#1)}#3}
\Crefrangeformat{equation}{Equations~\textup{#3(#1)#4--#5(#2)#6}}
\Crefmultiformat{equation}{Equations~\textup{#2(#1)#3}}{ and \textup{#2(#1)#3}}
{, \textup{#2(#1)#3}}{, and \textup{#2(#1)#3}}
\Crefrangemultiformat{equation}{Equations~\textup{#3(#1)#4--#5(#2)#6}}%
{ and \textup{#3(#1)#4--#5(#2)#6}}{, \textup{#3(#1)#4--#5(#2)#6}}{, and \textup{#3(#1)#4--#5(#2)#6}}
\crefdefaultlabelformat{#2\textup{#1}#3}

\crefname{chapter}{chapter}{chapters}
\crefname{appendix}{appendix}{appendices}
\crefname{subappendix}{section}{sections}
\Crefname{subappendix}{Section}{Sections}
\crefname{page}{page}{pages}
\crefname{subsection}{section}{sections}
\Crefname{subsection}{Section}{Sections}

\usepackage{csquotes}

\usepackage[
  mincrossrefs=99,
  style=numeric-comp,
  bibstyle=siam,
  sorting=nyt,
  backref=true,
  backend=biber,
  url=true,
  doi=true,
  maxbibnames=99,
  maxnames=99,
  giveninits=true,
  useprefix=true,
  date=year,
  autolang=other
]{biblatex}
\DeclareCiteCommand{\citeauthornsc}
  {%
   \boolfalse{citetracker}%
   \boolfalse{pagetracker}%
   \usebibmacro{prenote}}
  {\ifciteindex
     {\indexnames{labelname}}
     {}%
   \printnames{labelname}}
  {\multicitedelim}
  {\usebibmacro{postnote}}
\DeclareSourcemap{
  \maps[datatype=bibtex]{
    \map{
      \step[fieldsource=doi,final]
      \step[fieldset=url,null]
    }
  }
}

\usepackage{caption}
\usepackage{fancyhdr}
\usepackage{lastpage}

\title{Piecewise discretization of monodromy operators of delay equations on adapted meshes}

\author{Dimitri Breda, Davide Liessi and Rossana Vermiglio}

\newcommand{\mydate}{17 February 2022}

\fancypagestyle{mypagestyle}{%
  \fancyhf{}%
  \fancyhead[LO,RE]{\scriptsize Piecewise discretization of monodromy operators}%
  \fancyhead[RO,LE]{\scriptsize Breda, Liessi, Vermiglio}%
  \fancyfoot[LE,RO]{\scriptsize \thepage\,/\,\pagerefstar{LastPage}}%
  \fancyfoot[LO,RE]{\scriptsize \mydate{}}%
}
\begin{document}
\thispagestyle{empty}

\begin{center}
\LARGE
Piecewise discretization of monodromy operators of delay equations on adapted meshes

\bigskip
\large
Dimitri Breda\footnote{\email{dimitri.breda@uniud.it}},
Davide Liessi\footnote{\email{davide.liessi@uniud.it}},
Rossana Vermiglio\footnote{\email{rossana.vermiglio@uniud.it}}

\medskip
\small
CDLab -- Computational Dynamics Laboratory \\
Department of Mathematics, Computer Science and Physics, University of Udine \\
Via delle Scienze 206, 33100 Udine, Italy

\medskip
\large
\mydate
\end{center}

\begin{abstract}
Periodic solutions of delay equations are usually approximated as continuous piecewise polynomials on meshes adapted to the solutions' profile.
In practical computations this affects the regularity of the (coefficients of the) linearized system and, in turn, the effectiveness of assessing local stability by approximating the Floquet multipliers.
To overcome this problem when computing multipliers by collocation, the discretization grid should include the piecewise adapted mesh of the computed periodic solution.
By introducing a piecewise version of existing pseudospectral techniques, we explain why and show experimentally that this choice is essential in presence of either strong mesh adaptation or nontrivial multipliers whose eigenfunctions' profile is unrelated to that of the periodic solution.


\end{abstract}

\section{Introduction}

Periodic solutions and their asymptotic stability are among the prime interests in the study of dynamical systems.
In the case of delay differential equations (DDEs) these solutions are usually approximated with continuous piecewise polynomials determined by collocating a corresponding boundary value problem (BVP) on the period interval \cite{Bader1985,EngelborghsLuzyaninaIntHoutRoose2001}.
Then, as an accomplished standard (e.g., as in DDE-BIFTOOL%
\footnote{\url{http://ddebiftool.sourceforge.net/}}
\cite{EngelborghsLuzyaninaRoose2002,SieberEngelborghsLuzyaninaSamaeyRoose2014}), the partition of the period interval is adapted to the profile of the solution, moving away from uniform (for mesh adaptation see \cite{AscherMattheijRussell1988,EngelborghsLuzyaninaIntHoutRoose2001}).
Eventually, the local stability is assessed by computing the characteristic multipliers,%
\footnote{Computing stability indicators is essential in a continuation framework for bifurcation analysis.
Recomputing solutions after perturbation may represent an attractive alternative when parameters are fixed.}
relying on Floquet theory and on the principle of linearized stability, see \cite{DiekmannVanGilsVerduynLunelWalther1995,HaleVerduynLunel1993} for DDEs and \cite{BredaLiessi2021} for renewal equations (REs).

Recently, Borgioli et al.\ \cite{BorgioliHajduInspergerStepanMichiels2020} proposed a generalization of the collocation approach of DDE-BIFTOOL to compute Floquet multipliers of linear time-periodic DDEs, possibly with discontinuous coefficients.
They comment in section 4 that the points where the coefficients are not differentiable should be included in the collocation grid.
Accordingly, we highlight that when the linear system comes from linearizing a nonlinear problem around a numerically computed periodic solution, the resulting coefficients are in general only continuous, even for smooth problems, being the approximated solution a continuous piecewise polynomial.
This may deteriorate the convergence of the computated multipliers, as it happens, e.g., for the pseudospectral collocation methods for DDEs \cite{BredaMasetVermiglio2012}, for REs \cite{BredaLiessi2018} and for coupled REs and DDEs \cite{BredaLiessi2020}.
All these methods construct a matrix discretizing the monodromy operator by using a single polynomial on the whole domain interval.

Therefore, in this work we first recast the cited methods in a piecewise fashion, discretizing the monodromy operator on a grid including the adapted partition of the period interval from the given numerical periodic solution.
Then, with reference to the convergence analysis in \cite{BredaLiessi2018,BredaLiessi2020,BredaMasetVermiglio2012,BredaMasetVermiglio2015}, we explain why this choice is not only necessary to prevent order reduction or even loss of convergence, but also computationally convenient in the case of strongly adapted partitions.
Moreover, we discuss also the case of poor approximation of nontrivial multipliers whose eigenfunctions have large oscillations unrelated to the profile of the computed periodic solution \cite{YanchukRuschelSieberWolfrum2019}.
This leads to an increase of the error constants and a denser collocation grid is thus required.
In this respect, we show experimentally that this new grid should always be a refinement of the adapted mesh of the periodic solution: uniform grids with the same amount of nodes may fail to reach the desired accuracy.
The method we obtain applies to DDEs, REs and coupled REs and DDEs, with both discrete and distributed constant delays.

In this paper, after a summary in \cref{s_background} of the theoretical and numerical aspects related to approximating periodic solutions and studying their stability by computing the Floquet multipliers, we show in \cref{s_expected} an example of the difficulties encountered by the non-piecewise method \cite{BredaMasetVermiglio2012}.
Then, in \cref{s_reasonable}, we illustrate the piecewise reformulation, discuss its convergence, show with a simple DDE that order reduction may occur if the adapted partition is not taken into account and eventually repair the failure described in the previous \namecref{s_expected}.
In \cref{s_experimental} we provide other numerical experiments confirming the expected convergence behavior, the better performance in the case of solutions on strongly adapted partitions of the period interval, the suitability for multipliers relevant to eigenfunctions with smooth yet large oscillations and, finally, the versatility with respect to the classes of delay equations.

MATLAB codes implementing the described method are available at \url{http://cdlab.uniud.it/software}.

\section{Background}
\label{s_background}

In the following we summarize Floquet theory and local stability (\cref{s_floquet}), the numerical computation of periodic solutions (\cref{s_numericalper}) and the computation of the Floquet multipliers (\cref{s_numericalflo}), also introducing the necessary notations.
We restrict to DDEs and give references for REs and coupled equations.

\subsection{Floquet theory and local stability of periodic solutions}
\label{s_floquet}

Let $d_{Y}$ be a positive integer, $\tau$ a positive real and $\lvert\cdot\rvert$ any norm in finite dimension.
We consider DDEs
\begin{equation}\label{nonlineardde}
y'(t) = G(y_{t})
\end{equation}
for $G\colon Y \to \mathbb{R}^{d_{Y}}$, $Y \coloneqq C([-\tau, 0], \mathbb{R}^{d_{Y}})$ with norm $\lVert\psi\rVert_{Y} \coloneqq \max_{\theta \in [-\tau, 0]} \lvert\psi(\theta)\rvert$, and $y_{t}(\theta) \coloneqq y(t+\theta)$ for $\theta \in [-\tau, 0]$.

Assume that \cref{nonlineardde} has an $\omega$-periodic solution $\bar{y}$.
Linearizing \cref{nonlineardde} around $\bar{y}$ leads to the linear $\omega$-periodic DDE
\begin{equation}\label{linearizeddde}
y'(t) = DG(\bar{y}_{t}) y_{t}
\end{equation}
for $DG$ the Fréchet derivative of $G$.
Let $U(t,s) \colon Y\to Y$, $t\geq s$, be the associated evolution operator, i.e.,
\begin{equation*}
U(t, s) \psi = y(\cdot; s, \psi)_{t},
\end{equation*}
where $y(\cdot; s, \psi)$ is the solution of the initial value problem (IVP) for \cref{linearizeddde} with $y_{s} = \psi$ (the IVP is well-posed in the periodic case, see, e.g., \cite[Theorems 2.2.1, 2.2.2 and 2.2.3]{HaleVerduynLunel1993} and \cite[Theroem 3.7 and Remark 3.8]{Smith2011}).
The Floquet multipliers (simply multipliers in the sequel) are the eigenvalues of the monodromy operators $U(t+\omega,t)$, and they are independent of $t$ \cite[Theorem XIII.3.3]{DiekmannVanGilsVerduynLunelWalther1995}.
Note that $1$ is always a multiplier (usually called \emph{trivial}), since \cref{linearizeddde} is the linearization of a DDE around a periodic solution \cite[Theorem XIV.2.6]{DiekmannVanGilsVerduynLunelWalther1995}.
As is well known, the multipliers can give information on the local stability of $\bar{y}$ through the principle of linearized stability.
Namely, if $G$ is a $C^{1}$ function and the trivial multiplier is simple, then $\bar{y}$ is asymptotically stable if all the nontrivial multipliers are inside the unit circle; on the contrary, if there exists a nontrivial multiplier outside the unit circle, then $\bar{y}$ is unstable (for a proof follow \cite[Theorems XIV.3.3 and XIV.4.5]{DiekmannVanGilsVerduynLunelWalther1995}, with their hypotheses satisfied thanks to \cite[Exercise XIII.2.3 and section XIV.3]{DiekmannVanGilsVerduynLunelWalther1995}).

Similar results hold also for REs
\begin{equation*}
x(t) = F(x_{t})
\end{equation*}
with $F\colon X \to \mathbb{R}^{d_{X}}$, $X = L^{1}([-\tau, 0], \mathbb{R}^{d_{X}})$ with norm $\lVert \varphi\rVert_{X} \coloneqq \int_{-\tau}^{0} \lvert\varphi(\theta)\rvert \dd \theta$ and $d_{X}$ a positive integer.
Thanks to the Riesz representation theorem for $L^{1}$ (see, e.g., \cite[p. 400]{RoydenFitzpatrick2010}), the linearization around a possible periodic solution $\bar{x}$ has the form
\begin{equation*}
x(t) = \int_{-\tau}^{0} C(t, \theta) x_{t}(\theta) \dd \theta,
\end{equation*}
with $C\colon \mathbb{R} \times [-\tau, 0] \to \mathbb{R}^{d_{X} \times d_{X}}$ a measurable function, periodic in $t$.
Monodromy operators and relevant multipliers are defined as in the case of DDEs and the principle of linearized stability holds unchanged under mild regularity assumptions on $F$ (viz.\ $F$ is $C^{1}$ and globally Lipschitz continuous; the proof follows \cite[Theorems XIV.3.3 and XIV.4.5]{DiekmannVanGilsVerduynLunelWalther1995}, with their hypotheses satisfied thanks to \cite{BredaLiessi2018,BredaLiessi2021}).

Finally, for coupled equations
\begin{equation}\label{nonlinearcoupled}
\left\{
\begin{aligned}
& x(t) = F(x_{t}, y_{t}), \\
& y'(t) = G(x_{t}, y_{t}),
\end{aligned}
\right.
\end{equation}
where $F\colon X \times Y \to \mathbb{R}^{d_{X}}$, $G\colon X \times Y \to \mathbb{R}^{d_{Y}}$ and $X$ and $Y$ are as above with $\lVert(\varphi, \psi)\rVert_{X \times Y} \coloneqq \lVert\varphi\rVert_{X} + \lVert\psi\rVert_{Y}$, a Floquet theory is currently missing, but it is reasonable to expect the validity of similar results.

\subsection{Numerical computation of periodic solutions}
\label{s_numericalper}

In applications, exact periodic solutions of DDEs are in general unknown, so numerical methods are needed.
These typically consist in solving BVPs via piecewise orthogonal collocation \cite{Bader1985,EngelborghsLuzyaninaIntHoutRoose2001}.
Very recently, this methodology has been extended and applied to REs for the first time in \cite{BredaDiekmannLiessiScarabel2016}, but a systematic treatment for DDEs, REs and coupled equations appeared only in \cite{Ando2020}.
Concerning the convergence, the only available rigorous error analysis can be found in \cite{AndoBreda2020b} for DDEs and in \cite{Ando2021,AndoBreda} for REs.
Here we just summarize from \cite{EngelborghsLuzyaninaIntHoutRoose2001} the main aspects of this numerical scheme in the case of DDEs, which corresponds to the one implemented in DDE-BIFTOOL.
For the extension to REs and coupled equations see \cite{Ando2020}.

Assume again that \cref{nonlineardde} has a periodic solution $\bar{y}$, unknown together with its period $\omega$.
By rescaling the time through the map $s_{\omega}\colon\mathbb{R}\to\mathbb{R}$ defined as $s_{\omega}(t)\coloneqq t/\omega$, we can look at $\bar{y}$ as the solution of the BVP
\begin{equation}\label{bvp}
\left\{
\begin{aligned}
&y'(t)=\omega G(\tilde{y}_{t}\circ s_{\omega}),\quad t\in[0,1],\\
&y(0)=y(1),\\
&\phi(y)=0,
\end{aligned}
\right.
\end{equation}
where $\phi$ is a scalar (usually linear) function imposing a phase condition to remove translational invariance \cite{Doedel2007} and $\tilde{y}_{t}$ is defined as
\begin{equation*}
\tilde{y}_{t}(\theta)\coloneqq y(t+\theta+k), \quad t+\theta\in[-k,-k+1], \; k \in \mathbb{N},
\end{equation*}
exploiting the periodicity to evaluate the solution in $[-1,0]$ as required by the presence of the delay.

Let $L$ and $m$ be positive integers.
Consider a partition of $[0, 1]$ through $0 = t_{0} < t_{1} < \dots < t_{L} = 1$ and define the set of continuous piecewise $m$-degree polynomials
\begin{equation*}
\Pi_{L,m}\coloneqq\bigl\{p\in C([0,1],\mathbb{R}^{d_{Y}}) \mid p\restriction_{[t_{i},t_{i+1}]}\in\Pi_{m},\; i\in\{0,\dots,L-1\}\bigr\},
\end{equation*}
where $\Pi_{m}$ is the set of $\mathbb{R}^{d_{Y}}$-valued polynomials of degree at most $m$.
The piecewise collocation approach consists in looking for $p\in\Pi_{L,m}$ and $w\in\mathbb{R}$ satisfying
\begin{equation}\label{bvpcol1}
\left\{
\begin{aligned}
&p'(\zeta_{i,j})=wG(\tilde{p}_{\zeta_{i,j}}\circ s_{w}),\quad j\in\{1,\dots,m\},\;i\in\{0,\dots,L-1\},\\
&p(0)=p(1),\\
&\phi(p)=0,
\end{aligned}
\right.
\end{equation}
for a choice of $m$ collocation nodes $\zeta_{i,j}$ per interval, with $t_{i} \leq \zeta_{i, 1} < \dots < \zeta_{i, m} \leq t_{i+1}$, $i \in \{0, \dots, L-1\}$, typically the zeros of some family of orthogonal polynomials (e.g., Gauss--Legendre or Chebyshev).

As for the convergence of the method, it is shown in \cite{AndoBreda2020b} that, for $G$ sufficiently smooth, the collocation error $p-\bar{y}$ in the space of bounded measurable functions vanishes with order $m$ for $L\to\infty$ for the finite elements method (FEM, the most commonly adopted).
For the spectral elements method (SEM), i.e., for $m\to\infty$ and fixed $L$, experimental evidence of spectral accuracy \cite{Trefethen2000} is reported in \cite{Ando2020}, but no proof is available.

With regards to the implementation, \cref{bvpcol1} is recast as a system of nonlinear equations by using a suitable representation of the collocation polynomial.
The standard choice is the Lagrange form, defined by the basis of Lagrange polynomials $\ell_{i,0},\dots,\ell_{i,m}$ at the equidistant nodes $z_{i,j}\coloneqq t_{i}+j h_{i} / m$, for $i \in \{0,\dots,L-1\}$, $j\in\{0,\dots,m\}$ and $h_{i}:=t_{i+1}-t_{i}$.
The resulting system is typically solved by resorting to Newton's method, with a favorable arrow-shaped structure of the resulting Jacobian matrix (see \cite{Doedel2007} for a clear exposition in the case of ODEs).

As a final essential remark, we highlight again that standard implementations (as in DDE-BIFTOOL) make use of mesh adaptation: starting from a uniform partition of $[0,1]$, the distribution of the $L$ intervals is adapted to the solution's profile in order to control the overall error \cite{EngelborghsLuzyaninaIntHoutRoose2001}.
As we will base our collocation approach on the adapted mesh resulted from approximating $\bar{y}$, we introduce the ratio $\rho$ between the lengths of the largest and smallest intervals in the adapted partition as an indicator of how far it is from uniform.%
\footnote{Observe that $\rho\geq1$ and $\rho=1$ if and only if the partition is uniform.}
We anticipate that our approach does not perform any further adaptation on the solution's adapted mesh that it receives in input.

\subsection{Numerical computation of the Floquet multipliers}
\label{s_numericalflo}

The approach to computing the multipliers presented in \cite{BredaMasetVermiglio2012,BredaMasetVermiglio2015} for DDEs, \cite{BredaLiessi2018} for REs and \cite{BredaLiessi2020} for coupled equations is based on the non-piecewise discretization of a generic evolution operator following the relevant IVP, and it can thus be applied to approximate the spectrum of any such operator.
In this sense, this approach is more general than those used, e.g., in DDE-BIFTOOL or in \cite{BorgioliHajduInspergerStepanMichiels2020}, which explicitly exploit the structure of the periodic BVP (for other differences, see \cref{remark-B} below).
In this \namecref{s_numericalflo} we describe the discretization for DDEs presented in \cite{BredaMasetVermiglio2012}.
We refer the reader to the cited works and to \cite{Liessi2018} for the extension to REs and coupled equations.

Consider then \cref{linearizeddde} and the relevant evolution operator $T \coloneqq U(s+\omega, s)$ for $s \in \mathbb{R}$ and $\omega \geq 0$.
We apply pseudospectral collocation techniques to obtain a finite-dimensional approximation of $T$, by first conveniently reformulating it as follows.
Define the function spaces
$Y^{+} \coloneqq C([0, \omega], \mathbb{R}^{d_{Y}})$
and
$Y^{\pm} \coloneqq C([-\tau, \omega], \mathbb{R}^{d_{Y}})$
with the corresponding uniform norms.
Let $V\colon Y \times Y^{+} \to Y^{\pm}$ be the operator which, given an initial function $\psi$ on $[-\tau, 0]$ and the function $z$ prescribed by the right-hand side of \cref{linearizeddde} on $[0, \omega]$, constructs the solution of \cref{linearizeddde} on $[-\tau, \omega]$ as
\begin{equation*}
V(\psi,z)(t) \coloneqq
\begin{cases}
\displaystyle \psi(0) + \int_{0}^{t} z(\sigma) \dd\sigma, & t \in [0, \omega], \\
\psi(t), & t \in [-\tau, 0].
\end{cases}
\end{equation*}
Let also $\mathcal{F}_{s} \colon Y^{\pm} \to Y^{+}$ be the operator defined as
\begin{equation*}
(\mathcal{F}_{s}v) (t) \coloneqq DG(\bar{y}_{s+t})v_{t},\quad t\in[0,\omega],
\end{equation*}
which basically applies to its argument the action of the right-hand side of \cref{linearizeddde} (with the time shifted by $s$ so that the initial time is $0$).
Finally, $T$ can be reformulated as
\begin{equation}\label{T-as-V-dde}
T \psi = V(\psi,z^{\ast})_{\omega},
\end{equation}
where $z^{\ast}\in Y^{+}$ is the solution of the fixed point equation
\begin{equation}\label{fixed-point-dde}
z = \mathcal{F}_{s} V(\psi, z),
\end{equation}
which exists and is unique in the same conditions as the solutions of \cref{linearizeddde}.
Observe that $z^{\ast}$ is the derivative of the solution of \cref{linearizeddde} with initial function $y_{s}=\psi$.

\bigskip

Let now $M$ and $N$ be positive integers.
We consider partitions of $[- \tau, 0]$ and $[0, \omega]$ respectively defined by $- \tau = \theta_{M} < \dots < \theta_{0} = 0$ and $0 \leq t_{1} < \dots < t_{N} \leq \omega$.

The discretization of $Y$ is $Y_{M} \coloneqq \mathbb{R}^{d_{Y} (M + 1)}$, with elements
$\Psi = (\Psi_{0}, \dots, \Psi_{M})$ for $\Psi_{m} \in \mathbb{R}^{d_{Y}}$, $m \in \{0, \dots, M\}$.
We introduce the restriction operator $R_{M} \colon Y \to Y_{M}$ given by
$R_{M} \psi \coloneqq (\psi(\theta_{0}), \dots, \psi(\theta_{M}))$
and the prolongation operator $P_{M} \colon Y_{M} \to Y$ as the discrete Lagrange interpolation operator
$P_{M} \Psi(\theta) \coloneqq \sum_{m = 0}^{M} \ell_{m}(\theta) \Psi_{m}$, $\theta \in [-\tau, 0]$,
where $\ell_{0}, \dots, \ell_{M}$ are the Lagrange polynomials relevant to the nodes in $[-\tau, 0]$.
Observe that
\begin{equation}\label{PR-RP}
R_{M} P_{M} = I_{Y_{M}}, \quad\quad
P_{M} R_{M} = \mathcal{L}_{M},
\end{equation}
where $\mathcal{L}_{M} \colon Y \to Y$ is the Lagrange interpolation operator that associates to a function $\psi \in Y$ the $M$-degree $\mathbb{R}^{d_{Y}}$-valued polynomial $\mathcal{L}_{M} \psi$ such that
$\mathcal{L}_{M} \psi(\theta_{m}) = \psi(\theta_{m})$ for $m \in \{0, \dots, M\}$.
The discretization of $Y^{+}$ goes similarly by introducing $Y_{N}^{+}$, $R_{N}^{+}$, $P_{N}^{+}$ and $\mathcal{L}_{N}^{+}$ according to the partition of $[0, \omega]$.

Following \cref{T-as-V-dde,fixed-point-dde}, the discretization of $T$ is the finite-dimensional operator $T_{M, N} \colon Y_{M} \to Y_{M}$ defined as
\begin{equation}\label{discrete-T-as-V-dde}
T_{M, N} \Psi \coloneqq R_{M} V(P_{M} \Psi, P_{N}^{+} Z^{\ast})_{\omega},
\end{equation}
where $Z^{\ast} \in Y_{N}^{+}$ is the solution of the fixed point equation
\begin{equation}\label{discrete-fixed-point-dde}
Z = R_{N}^{+} \mathcal{F}_{s} V(P_{M}\Psi, P_{N}^{+}Z)
\end{equation}
for the given $\Psi \in Y_{M}$ (for the well-posedness see \cite{BredaMasetVermiglio2012}).
The eigenvalues of $T_{M, N}$ are then computed with standard methods and considered as approximations of the multipliers.

The convergence of the approximated multipliers has been proved in the cited works, and it holds under mild regularity assumptions on the ranges of $V$ and $\mathcal{F}_{s}$ (see, e.g., \cite[Theorem 3.3 and Proposition 4.5]{BredaMasetVermiglio2012}).
We summarize the main aspects in the following.
Let $\mu \in \mathbb{C} \setminus \{0\}$ be an eigenvalue of $T$ with generalized eigenspace $\mathcal{E}$, finite algebraic multiplicity $\nu$ and ascent $l$.
Let $\Delta$ be a neighborhood of $\mu$ such that $\mu$ is the only eigenvalue of $T$ in $\Delta$.
Then there exists a positive integer $\overline{N}$ such that, for any $N \geq \overline{N}$ and any $M \geq N$, $T_{M, N}$ has in $\Delta$ exactly $\nu$ eigenvalues $\mu_{M, N, j}$ for $j \in \{1, \dots, \nu\}$.
Moreover, if for each $\psi \in \mathcal{E}$ the function $z^{\ast}$ that solves \cref{fixed-point-dde} is of class $C^{p}$ for some $p \geq 1$, then
\begin{equation*}
\max_{j \in \{1, \dots, \nu\}} \lvert \mu_{M, N, j} - \mu\rvert = o\bigl(N^{\frac{1 - p}{l}}\bigr).
\end{equation*}
The result states that the order of convergence of the multipliers depends on the smoothness of the (derivative $z^{*}$ of the) solution of the IVP for \cref{linearizeddde} exiting from an eigenfunction of $\mu$.
Indeed, the approach collocates exactly this function $z^{\ast}$ on the period interval, and thus the error basically depends on the relvant interpolation procedure through the operators $R_{N}^{+}$, $P_{N}^{+}$ and $\mathcal{L}_{N}^{+}$.
This is the crucial point, and to simplify the following explanation let us assume that $G$ in \cref{nonlineardde} is smooth, so that also $\bar{y}$ is smooth (as periodicity combined with the smoothing effect of DDEs cancels possible breaking points).
Then, if we linearize \cref{nonlineardde} around the exact $\bar{y}$, it turns out that the coefficients of \cref{linearizeddde} are smooth and, as a consequence, so are the concerned eigenfunctions and the relevant functions $z^{\ast}$.%
\footnote{About the regularity of the eigenfunctions see respectively \cite{BredaMasetVermiglio2012} for DDEs and \cite{BredaLiessi2021} for REs.}
Then the order of convergence is infinite, which is coherent with the expectation that pseudospectral methods exhibit spectral accuracy \cite{Trefethen2000}.
On the contrary, if we linearize around a numerically computed $\bar{y}$, the coefficients of \cref{linearizeddde} are only continuous (yet piecewise analytic), and this lack of global smoothness can deteriorate the convergence behavior as we show in the next \namecref{s_expected}.

\begin{remark}
\label{remark-B}
We observe that the approach we presented, in the case of DDEs, collocates the derivative of the solution of the IVP, i.e., the solution $z^{\ast}$ of \cref{fixed-point-dde}, while DDE-BIFTOOL and \cite{BorgioliHajduInspergerStepanMichiels2020} collocate the solution of the periodic BVP \cref{bvp}: this may lead to slight differences in the experimental results even if using the same collocation grids.
As a further difference, we observe that the approach in \cite{BorgioliHajduInspergerStepanMichiels2020} does not discretize the monodromy operator (if not explicitly required), but recovers the multipliers from a generalized eigenvalue problem.
We note that we can follow the same alternative as well by separating the right-hand sides of \cref{discrete-T-as-V-dde,discrete-fixed-point-dde} with respect to the variables $\Psi$ and $Z$.
More comments are given at the end of \cref{s_reasonable}.
\end{remark}

\section{An expected failure}
\label{s_expected}

\begin{figure}[htp]
\begin{center}
\includegraphics{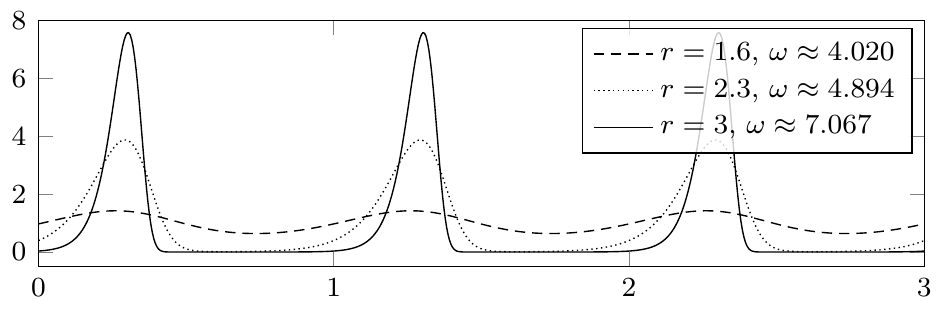}
\caption{Periodic solutions (rescaled to period $1$) of \cref{logisticdde} computed by DDE-BIFTOOL with $L_{\text{DB}} = 30$ and $m_{\text{DB}} = 6$.}
\label{fig:logisticddesol}
\end{center}
\end{figure}

To exemplify the approach described in \cref{s_background}, let us consider the delay logistic equation
\begin{equation}\label{logisticdde}
y'(t) = r y(t) (1-y(t-1)).
\end{equation}
To compute its periodic solutions we use DDE-BIFTOOL.%
\footnote{For all computations with DDE-BIFTOOL we use version 3.2a.}
Let $L_{\text{DB}}$ and $m_{\text{DB}}$ be, respectively, the number of pieces and the degree of the piecewise polynomials approximating such solutions.
We consider the solutions for $r \in \{1.6, 2.3, 3\}$ (\cref{fig:logisticddesol}) computed with $L_{\text{DB}} \in \{20, 30, 40\}$ and $m_{\text{DB}} \in \{2, 4, 6\}$.

In order to compute the relevant Floquet multipliers, we linearize \cref{logisticdde} around a generic solution $\bar{y}$, obtaining
\begin{equation*}
y'(t) = r (1-\bar{y}(t-1)) y(t) - r \bar{y}(t) y(t-1),
\end{equation*}
and apply the method of \cite{BredaMasetVermiglio2012}, implemented as \texttt{eigTMN} in the codes%
\footnote{\url{http://cdlab.uniud.it/software}}
accompanying \cite{BredaMasetVermiglio2015}, varying $M = N$.
We measure the error on both the trivial multiplier $1$ and the dominant nontrivial multiplier ($\mu \approx 0.8972$ for $r = 1.6$, $\mu \approx 0.1831 \cdot 10^{-2}$ for $r = 2.3$ and $\mu \approx 0.8037 \cdot 10^{-16} \pm 0.1198 \cdot 10^{-15} \mathrm{i}$ for $r = 3$).

\begin{figure}[htp]
\begin{center}
\includegraphics{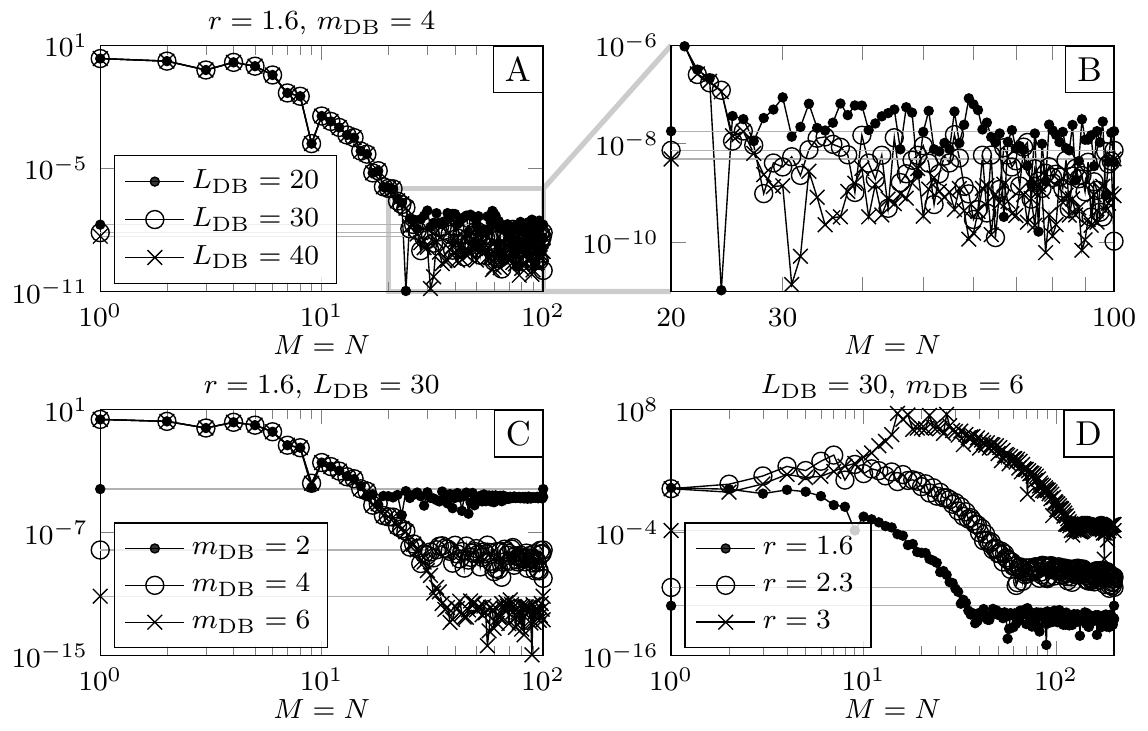}

\includegraphics{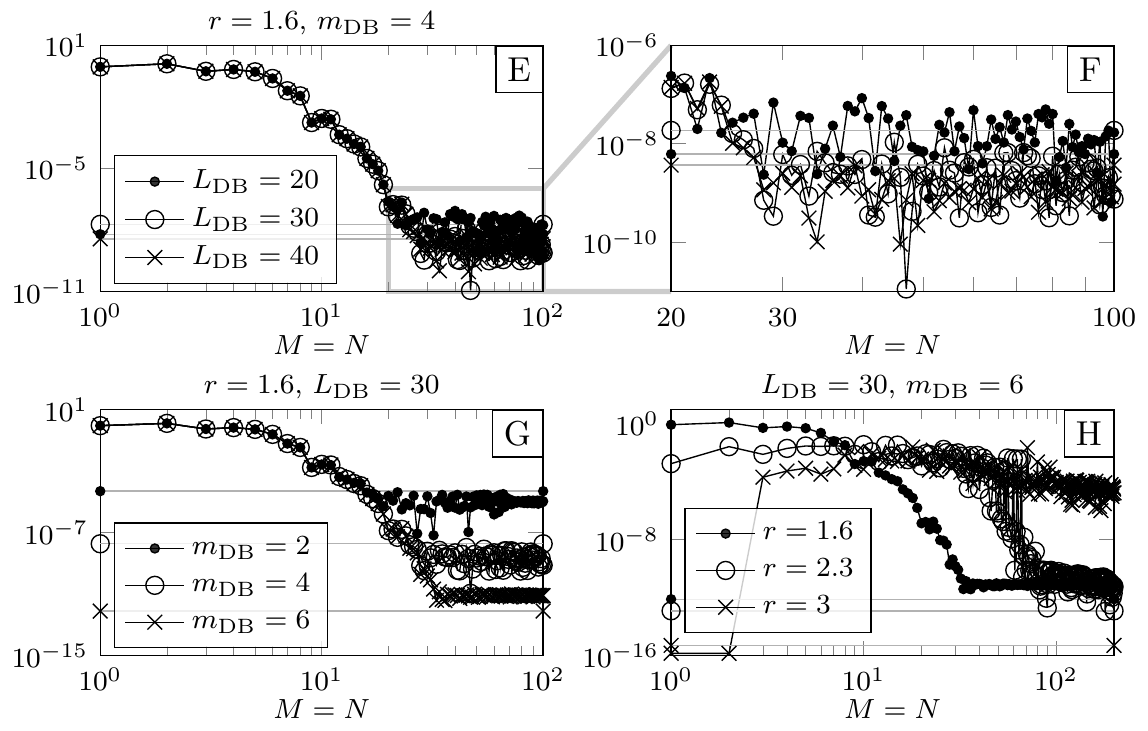}

\caption{Absolute errors of \texttt{eigTMN}, varying $M = N$, on the trivial (A--D) and dominant nontrivial (E--H) multipliers of \cref{logisticdde} linearized around the solutions computed by DDE-BIFTOOL.
The gray horizontal lines show DDE-BIFTOOL's errors on the same multipliers.
The reference values for the nontrivial multipliers are computed by DDE-BIFTOOL with $L_{\text{DB}} = 60$ and $m_{\text{DB}} = 10$.
\texttt{eigTMN}'s errors eventually decay with infinite order with barriers comparable to DDE-BIFTOOL's errors, except for the nontrivial multiplier with $r = 3$ (H), which seems to need even higher $M=N$.
For the trivial multiplier, as $r$ increases (D), the errors initially rise exponentially (exceeding $10^7$ for $r = 3$), complicating the choice of $M = N$.}
\label{fig:logisticddeerrboth}
\end{center}
\end{figure}

\Cref{fig:logisticddeerrboth} shows \texttt{eigTMN}'s errors as functions of $M = N$ varying $L_{\text{DB}}$, $m_{\text{DB}}$ and $r$ individually while keeping the others constant; for reference, they depict also DDE-BIFTOOL's errors on the same multipliers.
For the trivial multiplier (panels A--D) we eventually observe a convergence with infinite order with a barrier comparable to DDE-BIFTOOL's errors.
However, as $r$ increases (panel D), an initial phase of exponential rise appears, with errors exceeding $10^7$ in the worst case of $r = 3$.
This makes choosing an appropriate value for $M = N$ difficult.
As for the nontrivial multiplier (panels E--H), for $r = 1.3$ and $r = 2.6$ the convergence is similar, albeit lacking the initial exponential phase in the latter case, while for $r = 3$ \texttt{eigTMN} fails to reach DDE-BIFTOOL's accuracy with reasonable values of $M = N$.
In all cases we observe also that the influence of $L_{\text{DB}}$ on the error barrier is very limited, while the influence of $m_{\text{DB}}$ and $r$ is more significant.
For comparison, note that in this experiment \texttt{eigTMN}'s grid has up to $201$ nodes, while DDE-BIFTOOL's has $181$.

\begin{figure}[htp]
\begin{center}
\includegraphics{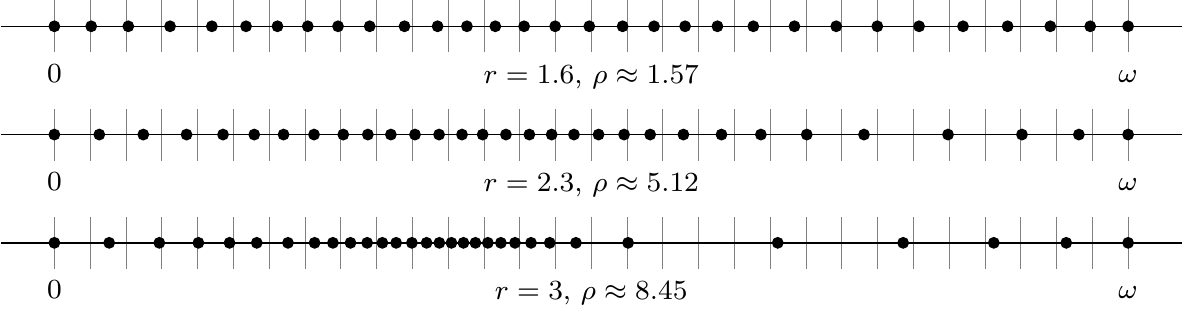}
\caption{Partitions of $[0, \omega]$ for the periodic solutions of \cref{fig:logisticddesol}, showing mesh adaptation as performed by DDE-BIFTOOL.
The vertical lines show the uniform partition.}
\label{fig:logisticddemeshes}
\end{center}
\end{figure}

Evidently, \texttt{eigTMN} has difficulties in dealing with this kind of problems.
We observe that as $r$ increases the solution (\cref{fig:logisticddesol}) transforms into a spike followed by a plateau, with larger and more sudden variations of the derivatives.
In such a case, if the solution's mesh is adapted, it progressively moves away from uniform.
As an example, \cref{fig:logisticddemeshes}
shows the partitions of the period interval for the solutions of \cref{fig:logisticddesol}.
The remedy we propose in \cref{s_reasonable} is to discretize the evolution operators on a collocation grid which includes the endpoints of the partition of the period interval on which the periodic solution has been computed.
In the following we will talk about more or less \emph{adapted} solutions and we will refer to the ratio $\rho$ introduced at the end of \cref{s_numericalper} to indicate how far from uniform their meshes have been adapted.

\section{A reasonable piecewise remedy}
\label{s_reasonable}

Let the partition of $[0, \omega]$ (as opposed to $[0, 1]$ in \cref{s_numericalper}) for the given numerical $\omega$-periodic solution be defined by $0 = t_{0} < t_{1} < \dots < t_{L} = \omega$.
As already noted in \cref{s_numericalper}, as collocation nodes the zeros of some family of orthogonal polynomials are usually chosen.
In our approach we collocate the operator also at the interval endpoints, which may need to be added to the collocation nodes.
Let them be $0 = c_{0} < \dots < c_{M} = 1$ in the interval $[0, 1]$; for each $i \in \{0, \dots, L-1\}$ and $j \in \{0, \dots, M\}$ we define $h_{i} \coloneqq t_{i+1} - t_{i}$ and $t_{i, j} \coloneqq t_{i} + h_{i} c_{j}$.
In our case, the collocation nodes typically are of Chebyshev type.

\begin{figure}[htp]
\begin{center}
\includegraphics{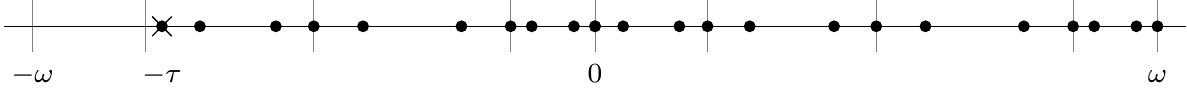}
\caption{Example collocation grid with $\omega > \tau$, $L=4$ and $M = 3$.
Ticks mark $t_{i}$ and $t_{i} - \omega$, the cross marks $-\tau$, dots mark the grid points.}
\label{fig:mesh}
\end{center}
\end{figure}

The interval $[-\tau, 0]$ too is partitioned according to the solution's mesh.
Assume that $\omega \geq \tau$.
The collocation nodes in $[-\tau, 0]$ are defined by subtracting $\omega$ to the nodes in $[0, \omega]$.
The leftmost resulting piece of $[-\tau, 0]$, however, requires special attention, since $-\tau$ may not coincide with one of the partition points.
Among the possible ways of treating it, the default choice of our implementation is to use $c_{0}, \dots, c_{M}$ to define new collocation nodes in the leftmost piece of $[-\tau, 0]$, independently of the corresponding nodes in $[0, \omega]$.
An example of the described collocation grid is depicted in \cref{fig:mesh}.
If instead $\omega < \tau$, the interval $[-\tau, 0]$ is partitioned in subintervals of length $\omega$ (with the leftmost possibly being smaller), each in turn partitioned according to the mesh of the numerical solution as described above.
For more details on the discretization, as mentioned in \cref{s_numericalflo}, we refer the reader to \cite{BredaLiessi2018,BredaMasetVermiglio2012} and, in particular for the piecewise approach, to \cite{BredaLiessiVermiglio}.

The discretization of the function spaces and of the operator $T$ is now almost identical to the one in \cref{s_numericalflo}, the only, but fundamental, difference being that the restriction, prolongation and Lagrange interpolation operators (now $R_{L,M}$, $P_{L,M}$ and $\mathcal{L}_{L,M}$, and $R_{L,M}^{+}$, $P_{L,M}^{+}$ and $\mathcal{L}_{L,M}^{+}$), act in a piecewise way.
We can thus obtain the discretization $T_{L,M}^{\text{pw}}$ of $T$ as
\begin{gather*}
T_{L,M}^{\text{pw}} \Psi \coloneqq R_{L,M} V(P_{L,M} \Psi, P_{L,M}^{+} Z^{\ast})_{\omega}, \\
Z = R_{L,M}^{+} \mathcal{F}_{s} V(P_{L,M}\Psi, P_{L,M}^{+}Z).
\end{gather*}
Recalling \cref{remark-B}, we note that the generalized eigenvalue problem possibly deriving from the piecewise version of the method results in sparse matrices, whose structure may be exploited for computational efficiency.

We implemented the method in the code \texttt{eigTMNpw}.%
\footnote{\url{http://cdlab.uniud.it/software}; for historical reasons we just added the suffix \texttt{pw} to the original name \texttt{eigTMN}.}
In the following we always use the default options; see \texttt{eigTMNpw}'s help for more details.
To avoid confusion with the discretization parameters of DDE-BIFTOOL ($L_{\text{DB}}$ and $m_{\text{DB}}$) and \texttt{eigTMN} ($M = N$), the number of pieces and the degree of the piecewise polynomials used in the collocation of the monodromy operator, respectively $L$ and $M$ in this \namecref{s_reasonable}, will be denoted as $L_{\text{pw}}$ and $M_{\text{pw}}$.
Moreover, where $L_{\text{pw}}$ is chosen to be the same as that of the given numerical periodic solution, it is intended that the solution's mesh in $[0, \omega]$ is used for \texttt{eigTMNpw} as well, unless otherwise indicated.

\subsection{About the convergence}
\label{s_about}

The differences in the formulation of the piecewise approach with respect to the non-piecewise approach of \cite{BredaLiessi2018,BredaLiessi2021,BredaMasetVermiglio2012,BredaMasetVermiglio2015,Liessi2018} are essentially limited to the restriction, prolongation and interpolation operators.
Most of the proofs of convergence in the cited works only depend on the essential properties \eqref{PR-RP} for the operators on $[0, \omega]$ ($R_{L,M}^{+}$, $P_{L,M}^{+}$ and $\mathcal{L}_{L,M}^{+}$), which are preserved in this new approach.
We thus expect the relevant convergence analysis to hold unchanged but for the underlying interpolation error that is at the basis of the convergence of the approximated multipliers as explained at the end of \cref{s_numericalflo} before \cref{remark-B}.
Indeed, as the interpolation process relies now on piecewise polynomials, we expect a convergence of finite order (proportional to the degree $M$ of the piecewise polynomials) when using the FEM and spectral accuracy \cite{Trefethen2000} when using the SEM.
Of course, this holds true only if the endpoints of the adapted partition of the period interval of the computed periodic solution are included in the collocation grid, so that interpolating polynomials are correctly used over pieces where the function to be interpolated is smooth.
Otherwise, as anticipated in \cref{s_numericalflo} and as shown experimentally in \cref{s_expected}, the overall convergence may be deteriorated.

\begin{figure}[htp]
\begin{center}
\includegraphics{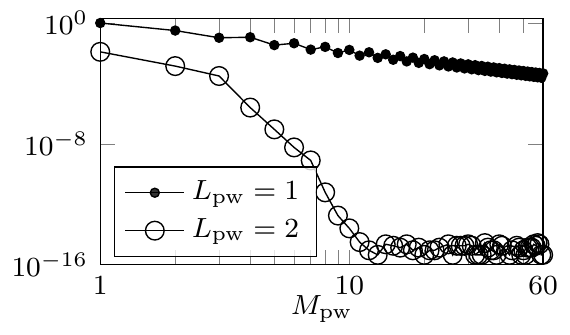}
\caption{Absolute errors of \texttt{eigTMNpw} on the dominant multiplier of \cref{tent}, whose coefficient is not differentiable at $1$.
The reference value is computed by \texttt{eigTMNpw} with $L_{\text{pw}} = 2$ and $M_{\text{pw}} = 120$.
When the mesh includes $1$ ($L_{\text{pw}} = 2$) the convergence order is infinite, otherwise ($L_{\text{pw}} = 1$) it is finite (precisely $2$).}
\label{fig:tent}
\end{center}
\end{figure}

This is illustrated by the next example, concerning the equation
\begin{equation}\label{tent}
x'(t) = (1-\lvert \mathop{\mathrm{mod}}(t,2) - 1\rvert) x(t-1),
\end{equation}
whose coefficient is piecewise linear, has period $2$ and is not differentiable at integer values of $t$.
\Cref{fig:tent} shows the errors on the dominant multiplier ($\mu \approx 2.0133$) when the partition of $[0, 2]$ does ($L_{\text{pw}}=2$) or does not ($L_{\text{pw}}=1$) contain $1$: indeed, in the former case the order of convergence is infinite, while in the latter it is finite.

Returning now to the example of \cref{s_expected}, in particular to panels D and H of \cref{fig:logisticddeerrboth}, we observe that the errors of \texttt{eigTMNpw} using the solution's mesh are comparable to those of DDE-BIFTOOL (\cref{tab:logisticddefixedreigtmnpw}).
\begin{table}[htp]
\begin{center}
\sisetup{
round-mode=figures,
round-precision=4,
scientific-notation=true,
table-format=1.3e+1,
}%
\begin{tabular}{S[scientific-notation=false,round-mode=off,table-format=1.1]SSSS}
\toprule
& \multicolumn{2}{c}{trivial multiplier} & \multicolumn{2}{c}{dominant nontrivial multiplier} \\
\cmidrule(rl){2-3}\cmidrule(rl){4-5}
$r$ & {DDE-BIFTOOL} & \texttt{eigTMNpw} & {DDE-BIFTOOL} & \texttt{eigTMNpw} \\
\midrule
1.6 & 7.26974036524553e-12 & 9.35251875944232e-13 & 7.88924481298636e-13 & 8.94251339644825e-12 \\
2.3 & 4.46265691067538e-10 & 2.44393172366131e-10 & 1.24292720213304e-13 & 6.59657202993491e-12 \\
3 & 0.000157671671576876 & 0.000344299146424198 & 5.08187323814298e-16 & 3.96047033460188e-15 \\
\bottomrule
\end{tabular}
\smallskip

\caption{Absolute errors of \texttt{eigTMNpw} and DDE-BIFTOOL on the trivial and dominant nontrivial multipliers of \cref{logisticdde} linearized around the solutions computed by DDE-BIFTOOL with $L_{\text{DB}} = 30$ and $m_{\text{DB}} = 6$.
\texttt{eigTMNpw} uses $L_{\text{pw}} = L_{\text{DB}}$ and $M_{\text{pw}} = m_{\text{DB}}$.
The reference values for the nontrivial multipliers are computed by DDE-BIFTOOL with $L_{\text{DB}} = 60$ and $m_{\text{DB}} = 10$.
The errors are comparable in magnitude.}
\label{tab:logisticddefixedreigtmnpw}
\end{center}
\end{table}
This is particularly notable for the nontrivial multiplier with $r = 3$, whose convergence in the non-piecewise case is very slow.
In all other cases we also note that to reach a comparable accuracy DDE-BIFTOOL and \texttt{eigTMNpw} use $181$ collocation nodes, while \texttt{eigTMN} needs less than $120$.
In fact, we show in \cref{s_stronglyadapted} that the piecewise approach becomes computationally convenient for solutions with higher values of $\rho$.

\bigskip

Finally, to close the discussion on convergence, the following observation should be taken into due consideration, also in view of tackling the issue of possible large oscillating eigenfunctions reported in \cite{YanchukRuschelSieberWolfrum2019}.
If including the endpoints of the adapted mesh of the computed periodic solution into the collocation grid is necessary to preserve the desired convergence order, still the error constants relevant to each piece depend on the derivatives of the interpolated function on that piece.
When this function is smooth but large derivatives appear, the error constants increase and a denser discretization should be considered.
However, for what explained above, the denser grid should always be a refinement of the adapted mesh from the periodic solution.
Numerical evidence of this fact will be given in \cref{s_nontrivial}.

\section{Experimental validation}
\label{s_experimental}

The experiments in \cref{s_testingconvergence} show the convergence properties of \texttt{eigTMNpw} by applying it to an equation with an explicitly known periodic solution, which is thus only affected by rounding errors.
\Cref{s_stronglyadapted} presents the application to an extreme case of strongly adapted solution, highlighting the advantage of the piecewise approach.
The case of eigenfunctions with oscillations unrelated to the profile of the periodic solution, anticipated in the introduction, is treated in \cref{s_nontrivial}.
Finally, in \cref{s_coupled} a coupled equation is considered, confirming the versatility of the proposed piecewise pseudospectral approach also in this piecewise reformulation.

\subsection{Testing the convergence}
\label{s_testingconvergence}

In order to show the convergence properties of \texttt{eigTMNpw}, we consider the RE with quadratic nonlinearity
\begin{equation}\label{quadraticre}
x(t) = \frac{\gamma}{2} \int_{-3}^{-1} x(t+\theta) (1-x(t+\theta)) \dd\theta,
\end{equation}
which has a branch of periodic solutions with the explicit expression
\begin{equation}\label{quadraticresolution}
\bar{x}(t) = \frac{1}{2} + \frac{\pi}{4 \gamma} + \sqrt{\frac{1}{2} - \frac{1}{\gamma} - \frac{\pi}{2 \gamma^{2}}\Bigl(1+\frac{\pi}{4}\Bigr)} \sin\Bigl(\frac{\pi}{2} t\Bigr),
\end{equation}
as proved in \cite{BredaDiekmannLiessiScarabel2016}.
To study the stability of $\bar{x}$, we consider the linear RE%
\footnote{\Cref{quadraticrelin} is not actually the linearization of \cref{quadraticre} in $L^1$, although it can be used to study the stability of its equilibria.
See \cite[section 3.5]{DiekmannGettoGyllenberg2008} for details; the extension of the results therein to periodic solutions is an open problem.}
\begin{equation}\label{quadraticrelin}
x(t) = \frac{\gamma}{2} \int_{-3}^{-1} (1-2\bar{x}(t+\theta)) x(t+\theta) \dd\theta.
\end{equation}
We use \texttt{eigTMNpw} to compute the multipliers relevant to \cref{quadraticresolution} for $\gamma = 4$ varying $L_{\text{pw}}$ (here defining a uniform partition of $[0, \omega]$) and $M_{\text{pw}}$.

\begin{figure}[htp]
\begin{center}
\includegraphics{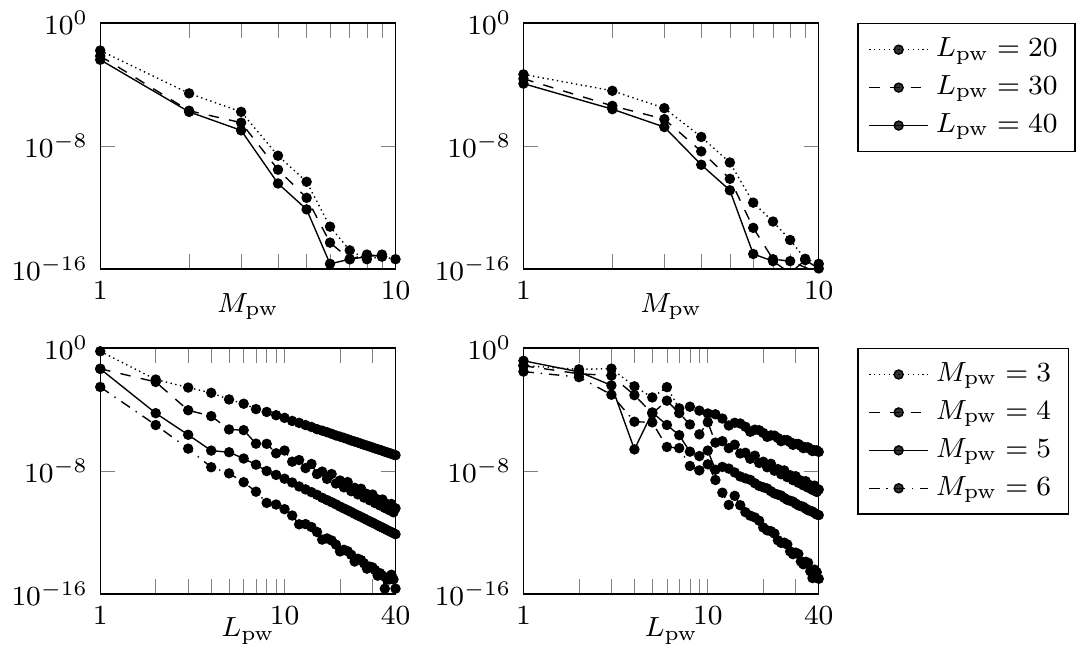}
\caption{Absolute errors of \texttt{eigTMNpw} on the trivial (left) and on the dominant nontrivial (right) multipliers of \cref{quadraticrelin} with $\gamma = 4$ in \cref{quadraticresolution}.
The reference value for the nontrivial multiplier is computed with $L_{\text{pw}} = 40$ and $M_{\text{pw}} = 15$.
The convergence order is infinite for the SEM (top) and finite for the FEM (bottom, precisely $4$, $6$, $6$ and $8$ on the left and $4$, $6$, $6$ and $10$ on the right).}
\label{fig:quadraticreerr}
\end{center}
\end{figure}

\Cref{fig:quadraticreerr} show the errors on the trivial multiplier (left panels) and on the dominant nontrivial multiplier ($\mu \approx -0.1355$, right panels), confirming our expectations (see \cref{s_about}): the error of the SEM (top panels) vanishes with infinite order, while the error of the FEM (bottom panels) vanishes with finite order (sometimes higher than the theoretical bound), increasing with $M_{\text{pw}}$.

\subsection{Strongly adapted solutions}
\label{s_stronglyadapted}

\begin{figure}[htp]
\begin{center}
\includegraphics{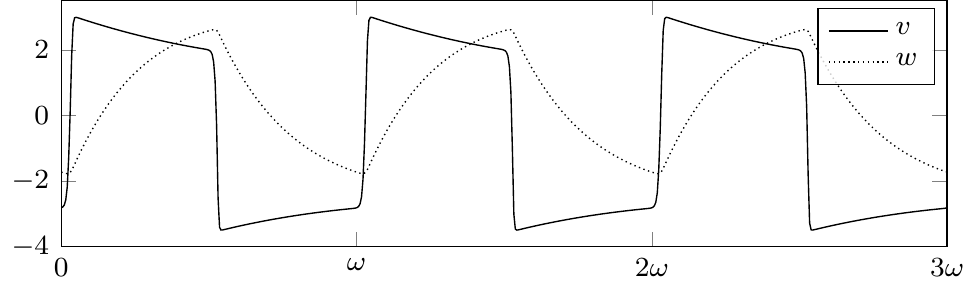}
\caption{Periodic solution of \cref{plantneural} with $a = 0.7$, $b = 0.8$, $\eta = -2$, $r = 0.08$ and $\tau = 25$ ($\omega \approx 50.7326$, $v_{0} \approx -1.1994$), computed by DDE-BIFTOOL with $L_{\text{DB}} = 30$ and $m_{\text{DB}} = 5$.}
\label{fig:plantneuralsol}
\end{center}
\end{figure}

\begin{figure}[htp]
\begin{center}
\includegraphics{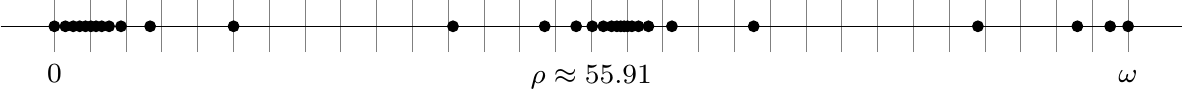}
\caption{Partition of $[0, \omega]$ for the solution of \cref{fig:plantneuralsol}.
The vertical lines show the uniform partition.}
\label{fig:plantneuralmesh}
\end{center}
\end{figure}

We consider the DDE
\begin{equation}\label{plantneural}
\left\{
\begin{aligned}
& v'(t) = v(t) - \frac{v(t)^{3}}{3} - w(t) + \eta (v(t-\tau) - v_{0}), \\
& w'(t) = r (v(t) + a - b w(t)), \\
& v_{0} \text{ a real root of } v - \frac{v^3}{3} - \frac{v+a}{b},
\end{aligned}
\right.
\end{equation}
proposed by Plant in \cite{Plant1981} to model recurrent neural feedback; its periodic solutions were studied in \cite{CastelfrancoStech1987}.
Note that $v_{0}$ is unique if $a \neq 0$ and $0 < b \leq 1$.
\Cref{fig:plantneuralsol} shows a periodic solution of \cref{plantneural} computed by DDE-BIFTOOL on a rather extremely adapted mesh (\cref{fig:plantneuralmesh}, $\rho \approx 55.91$): indeed it was used in \cite{EngelborghsLuzyaninaIntHoutRoose2001} to demonstrate the collocation method for computing periodic solutions with adaptive mesh selection which became part of DDE-BIFTOOL's foundations.

We linearize \cref{plantneural} around a generic solution $(\bar{v}, \bar{w})$, obtaining
\begin{equation*}
\left\{
\begin{aligned}
& v'(t) = (1 - \bar{v}(t)^{2}) v(t) - w(t) + \eta v(t-\tau), \\
& w'(t) = r v(t) - r b w(t).
\end{aligned}
\right.
\end{equation*}
We compute the solution with DDE-BIFTOOL for different values of $L_{\text{DB}}$ and $m_{\text{DB}}$.
For each solution we compute the relevant Floquet multipliers with \texttt{eigTMN} for increasing $M = N$ and with \texttt{eigTMNpw} with $L_{\text{pw}} = L_{\text{DB}}$ and $M_{\text{pw}} = m_{\text{DB}}$.

\begin{figure}[htp]
\begin{center}
\includegraphics{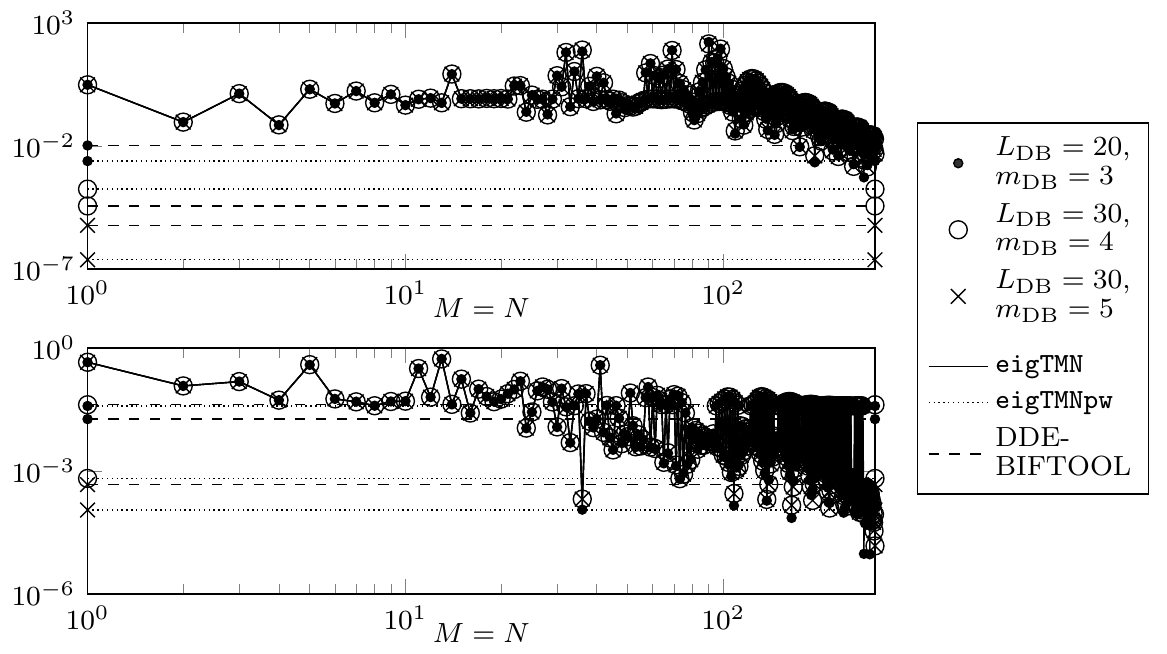}
\caption{Absolute errors of \texttt{eigTMN} for varying $M = N$ on the trivial (top) and dominant nontrivial (bottom) multipliers of \cref{plantneural} linearized around the solutions computed by DDE-BIFTOOL (parameters as in \cref{fig:plantneuralsol}), compared to the errors on the same multipliers of \texttt{eigTMNpw} with $L_{\text{pw}} = L_{\text{DB}}$ and $M_{\text{pw}} = m_{\text{DB}}$ and, for reference, of DDE-BIFTOOL (horizontal lines).
The reference value for the nontrivial multiplier is computed by DDE-BIFTOOL with $L_{\text{DB}} = 60$ and $m_{\text{DB}} = 10$.
While \texttt{eigTMNpw}'s errors are comparable with those of DDE-BIFTOOL, \texttt{eigTMN}'s ones either require very large $M = N$ to reach the desired magnitude (trivial multiplier) or exhibit an alternating behavior which makes choosing an appropriate $M = N$ difficult (nontrivial multiplier).}
\label{fig:plantneuralerrboth}
\end{center}
\end{figure}

\Cref{fig:plantneuralerrboth} shows the corresponding errors on the trivial multiplier (top panel) and on the dominant nontrivial multiplier ($\mu \approx 0.1444 \pm 0.0382\mathrm{i}$, bottom panel) compared with DDE-BIFTOOL's errors.
For the trivial multiplier, even with $M = N = 300$ the errors of \texttt{eigTMN} barely reach the same order of magnitude as DDE-BIFTOOL's \emph{worst} error; for the nontrivial multiplier, the errors alternate between largely different magnitudes, making it difficult to choose an appropriate $M = N$.
On the other hand, the errors of \texttt{eigTMNpw} are comparable to the ones of DDE-BIFTOOL, or even better in some cases.
The piecewise approach uses $151$ nodes and the non-piecewise approach uses up to $301$ nodes, which suggests that in this case the piecewise approach is also computationally more convenient.

\bigskip

\begin{figure}[htp]
\begin{center}
\includegraphics{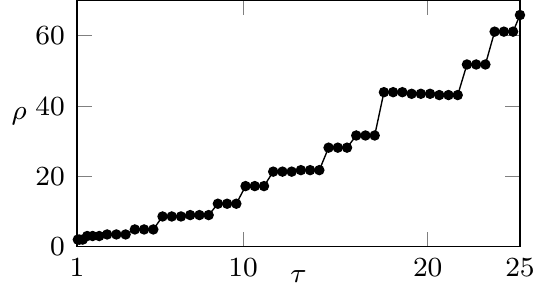}
\caption{Value of the ratio $\rho$ for the solutions of \cref{plantneural} computed by DDE-BIFTOOL with $L_{\text{DB}} = 30$ and $m_{\text{DB}} = 5$ for varying $\tau$ (other parameters as in \cref{fig:plantneuralsol}).}
\label{fig:plantneuralallratio}
\end{center}
\end{figure}

\begin{figure}[htp]
\begin{center}
\includegraphics{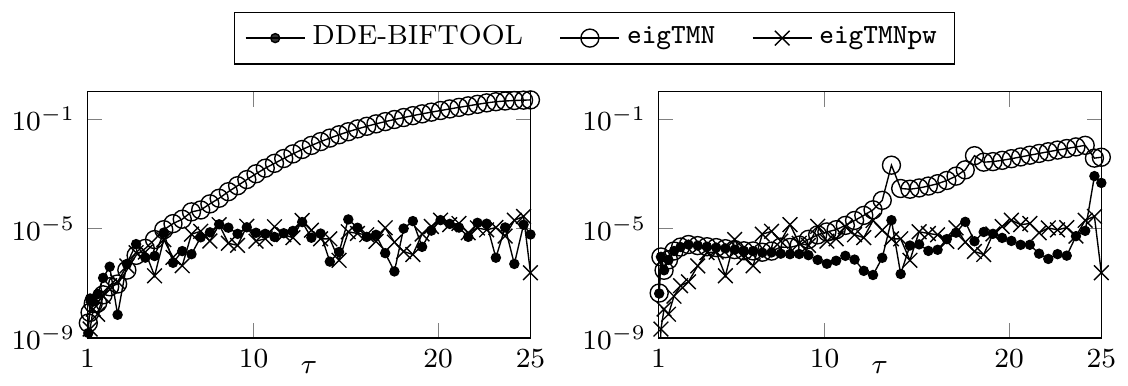}
\caption{Absolute errors of \texttt{eigTMN} and \texttt{eigTMNpw}, compared to those of DDE-BIFTOOL, on the trivial (left) and dominant nontrivial (right) multipliers of \cref{plantneural} linearized around the solutions computed by DDE-BIFTOOL with $L_{\text{DB}} = 30$, $m_{\text{DB}} = 5$ and varying $\tau$ (other parameters as in \cref{fig:plantneuralsol}).
In all cases $L_{\text{pw}} = L_{\text{DB}}$ and $M_{\text{pw}} = m_{\text{DB}}$ are used for \texttt{eigTMNpw}, while the errors of \texttt{eigTMN} are the mean error for $M = N \in \{140, \dots, 160\}$.
The reference value for the nontrivial multiplier is computed by DDE-BIFTOOL with $L_{\text{DB}} = 60$ and $m_{\text{DB}} = 10$.
\texttt{eigTMNpw}'s errors behave similarly to those of DDE-BIFTOOL, while \texttt{eigTMN}'s ones gradually increase.}
\label{fig:plantneuralallerr}
\end{center}
\end{figure}

In \cref{s_about} we anticipated that the computational convenience of the piecewise approach depends on the value of $\rho$.
The next experiment shows that this is indeed the case: we compute several solutions of \cref{plantneural} for varying $\tau$ and compare the errors on the multipliers as computed by DDE-BIFTOOL and \texttt{eigTMNpw} with $L_{\text{pw}} = L_{\text{DB}} = 30$ and $M_{\text{pw}} = m_{\text{DB}} = 5$, i.e., using $151$ nodes, and by \texttt{eigTMN} with $M = N \in \{140, \dots, 160\}$ (we consider the mean error in this case).
We first observe in \cref{fig:plantneuralallratio} that $\rho$ is almost monotonically increasing as $\tau$ varies from $1$ to $25$; we can thus use $\tau$ as a proxy for $\rho$.
In \cref{fig:plantneuralallerr} we observe that as $\tau$ increases DDE-BIFTOOL and \texttt{eigTMNpw} are equally accurate, while \texttt{eigTMN} progressively loses accuracy.

\subsection{Nontrivial multipliers with oscillating eigenfunctions}
\label{s_nontrivial}

As anticipated in the introduction, an eigenfunction may present large oscillations unrelated to the profile of the periodic solution \cite{YanchukRuschelSieberWolfrum2019}.
This necessarily requires the use of a denser discretization grid as explained in \cref{s_about}.
A key aspect is that the denser grid needs to be a refinement of the piecewise mesh of the numerical periodic solution, thus including the endpoints of pieces, since at those points the coefficients of the linearized equation are not smooth.

\begin{figure}[htp]
\begin{center}
\includegraphics{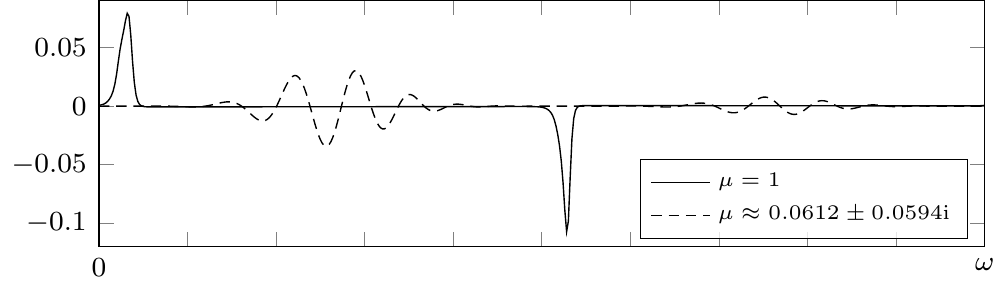}
\caption{Eigenfunctions relevant to the multiplier $\mu$ of \cref{plantneural} linearized around the solution computed by DDE-BIFTOOL with $L_{\text{DB}} = 240$ and $m_{\text{DB}} = 10$ (parameters as in \cref{fig:plantneuralsol}).}
\label{fig:plantneuralfinereigfun}
\end{center}
\end{figure}

\begin{figure}[ht]
\begin{center}
\includegraphics{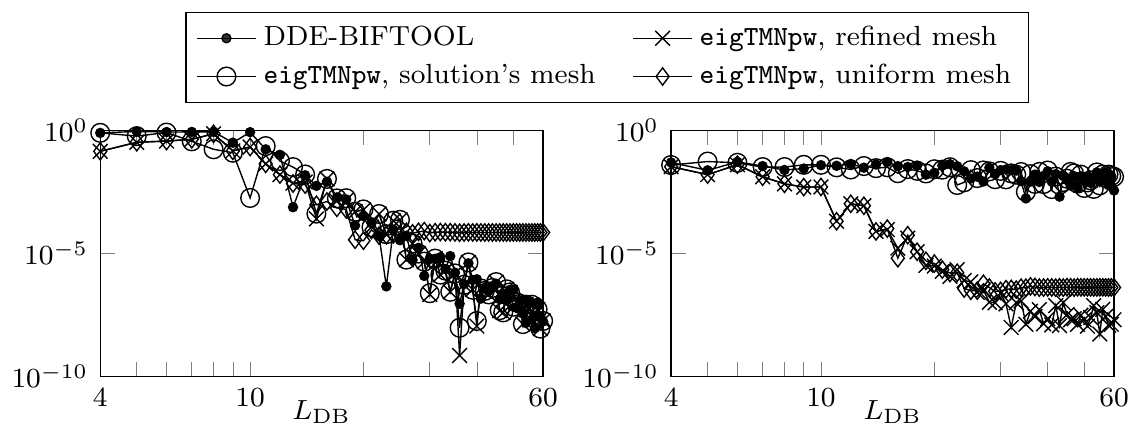}
\caption{Absolute errors on the trivial multiplier (left) and on the nontrivial multiplier $\mu \approx 0.0612 \pm 0.0594\mathrm{i}$ (right) of \cref{plantneural} linearized around the solutions computed by DDE-BIFTOOL with $m_{\text{DB}} = 5$ and varying $L_{\text{DB}}$ (parameters as in \cref{fig:plantneuralsol}).
For \texttt{eigTMNpw} $[0, \omega]$ is partitioned using the solution's mesh, a refinement of the latter (resulting in $129 \leq L_{\text{pw}} \leq 153$), and a uniform mesh with $L_{\text{pw}} = 153$ ($M_{\text{pw}} = m_{\text{DB}}$ in all cases).
DDE-BIFTOOL uses the solution's mesh.
The reference value for the nontrivial multiplier is computed by DDE-BIFTOOL with $L_{\text{DB}} = 240$ and $m_{\text{DB}} = 10$.
If the eigenfunction has oscillations unrelated to the solution's mesh (right), using the latter seems to prevent the convergence, while a denser partition allows to reach fairly small error barriers, smaller if it is actually a refinement.
For the trivial multiplier (left), the solution's mesh is well adapted also to the eigenfunction: the error vanishes similarly with both the original and the finer partition, while with the uniform one it reaches a barrier larger than in the other case.}
\label{fig:plantneuralfiner}
\end{center}
\end{figure}

To exemplify this fact, we turn our attention again to \cref{plantneural} and consider a nontrivial multiplier ($\mu \approx 0.0612 \pm 0.0594\mathrm{i}$) whose eigenfunction oscillates where the periodic solution and the eigenfunction of the trivial multiplier are almost constant and the adapted mesh has few points (see \cref{fig:plantneuralfinereigfun} and recall \cref{fig:plantneuralmesh}).
\Cref{fig:plantneuralfiner} compares the errors of DDE-BIFTOOL with those of \texttt{eigTMNpw} partitioning the period interval in three ways: using the solution's mesh, a refinement of the latter, and a dense but uniform mesh.%
\footnote{In the latter two cases the partitions are such that no piece is longer than five times the length of the longest piece of the solution's mesh.}
For the nontrivial multiplier (right panel), both denser partitions help in achieving the convergence of the multiplier, with the error barrier of the refined one being slightly better than that of the uniform one.
For the trivial multiplier (left panel), on the other hand, the uniform partition introduces an error barrier larger than the one for the nontrivial multiplier: this is expected, since the uniform partition is actually sparser than the adapted one where the periodic solution and the eigenfunction of the trivial multiplier have large oscillations.

\subsection{Coupled equations}
\label{s_coupled}

The last experiment shows that the pseudospectral approach in general, and the method we described in particular, are very versatile in terms of equation classes.
Indeed, in the previous \namecrefs{s_stronglyadapted} we considered an RE and a DDE, while in this \namecref{s_coupled} we consider a coupled equation.

In order to use an equation with a strongly adapted solution, we derive a coupled equation from \cref{plantneural} by integrating the equation for $w$, resulting in
\begin{equation}\label{plantneuralcoupled}
\left\{
\begin{aligned}
& v'(t) = v(t) - \frac{v(t)^{3}}{3} - w(t) + \eta (v(t-\tau) - v_{0}), \\
& w(t) = w(t-\tau) + \int_{t-\tau}^{t} r (v(s) + a - b w(s)) \dd s, \\
\end{aligned}
\right.
\end{equation}
where $v_{0}$ is the same as in \cref{plantneural}.

Note that, due to the presence in the RE for $w$ of the evaluation of $w$ at a specific point, \cref{plantneuralcoupled} does not belong to the family described by \cref{nonlinearcoupled} since the state space for $w$ cannot be $L^{1}$.
In fact \cref{plantneuralcoupled} is an example of neutral RE.
The method described in this work has been implemented in \texttt{eigTMNpw} with an eye to dealing also with neutral REs, but there is currently no proof of convergence for this case: the numerical treatment of neutral REs is the subject of ongoing research.
For more details on neutral REs and the relevant perturbation theory see the recent work \cite{DiekmannVerduynLunel2021}.

Linearizing \cref{plantneuralcoupled} around a generic solution $(\bar{v}, \bar{w})$, we obtain
\begin{equation*}
\left\{
\begin{aligned}
& v'(t) = (1 - \bar{v}(t)^{2}) v(t) - w(t) + \eta v(t-\tau), \\
& w(t) = w(t-\tau) + r \int_{t-\tau}^{t} v(s) \dd s - r b \int_{t-\tau}^{t} w(s) \dd s.
\end{aligned}
\right.
\end{equation*}
We consider the same periodic solution computed with DDE-BIFTOOL in \cref{s_stronglyadapted} with $L_{\text{DB}} = 30$ and $m_{\text{DB}} = 5$.
The relevant Floquet multipliers are computed using \texttt{eigTMNpw} first with $L_{\text{pw}} = 1$ for increasing $M_{\text{pw}}$ to show the failure of the non-piecewise approach, and then with $L_{\text{pw}} = L_{\text{DB}}$ and $M_{\text{pw}} = m_{\text{DB}}$.

\begin{figure}[htp]
\begin{center}
\includegraphics{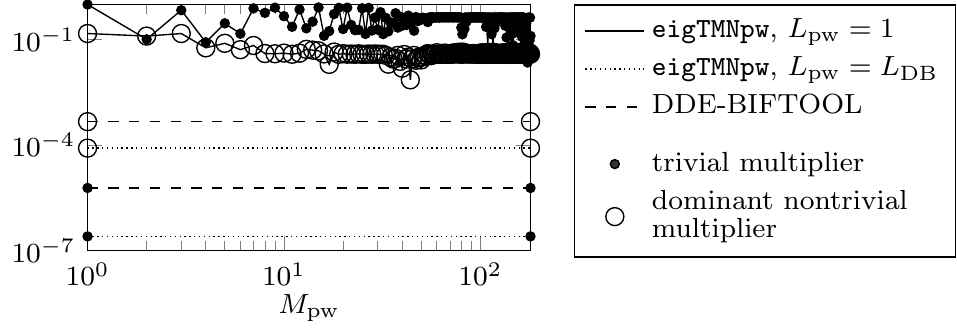}
\caption{Absolute errors of \texttt{eigTMNpw} on the trivial and dominant nontrivial multipliers of \cref{plantneuralcoupled} linearized around the solution of \cref{plantneural} computed by DDE-BIFTOOL with $L_{\text{DB}} = 30$ and $m_{\text{DB}} = 5$ (parameters as in \cref{fig:plantneuralsol}).
\texttt{eigTMNpw} is used in a nonpiecewise fashion for varying $M_{\text{pw}}$ and in a piecewise fashion with $L_{\text{pw}} = L_{\text{DB}}$ and $M_{\text{pw}} = m_{\text{DB}}$.
The reference value for the nontrivial multiplier is computed by DDE-BIFTOOL for \cref{plantneural} with $L_{\text{DB}} = 60$ and $m_{\text{DB}} = 10$.
\texttt{eigTMNpw}'s errors are compared to those of DDE-BIFTOOL for \cref{plantneural}: in the piecewise case the former are even more accurate than the latter.}
\label{fig:plantneuralcouplederr}
\end{center}
\end{figure}

\Cref{fig:plantneuralcouplederr} shows the corresponding errors on the trivial and dominant nontrivial multipliers compared with DDE-BIFTOOL's errors.
We observe that the errors of the non-piecewise approach suggest that the multipliers may begin to converge for high values of $M_{\text{pw}}$, with the error being greater than $10^{-2}$ for all values of $M_{\text{pw}} \leq 200$ (except one).
On the other hand, the error of the piecewise approach is smaller than DDE-BIFTOOL's error (in this case $151$ nodes are used).

\begin{remark}
In \cite{LuzyaninaEngelborghs2002} Luzyanina and Engelborghs compute the periodic solutions and the corresponding multipliers with DDE-BIFTOOL and experimentally study the convergence of the FEM.
In some of their examples the trivial multiplier converges with infinite order, while the order is finite for the nontrivial ones.
We note here that in our experience with \texttt{eigTMN} and \texttt{eigTMNpw} we never observed the superconvergence of the trivial multiplier.
\end{remark}

\section*{Acknowledgments}

The authors are members of INdAM Research group GNCS and of UMI Research group ``Mo\-del\-li\-sti\-ca socio-epidemiologica''.
This work was partially supported by the Italian Ministry of University and Research (MUR) through the PRIN 2020 project (No.\ 2020JLWP23) ``Integrated Mathematical Approaches to Socio-Epidemiological Dynamics'' (CUP: E15F21005420006).
The work of Davide Liessi was partially supported by Finanziamento Giovani Ricercatori 2020--2021 of INdAM Research group GNCS.

{\sloppy
\printbibliography
\par}

\end{document}